\documentclass{gtart}

\usepackage{amsmath,amssymb,latexsym,rlepsf}
\let\relabela\adjustrelabel

\input gtoutput
\volumenumber{3}\papernumber{12}\volumeyear{1999}
\pagenumbers{269}{302}\published{11 September 1999}
\proposed{Walter Neumann}
\seconded{Wolfgang Metzler, Joan Birman}
\received{27 November 1998}\revised{5 August 1999}
\accepted{5 September 1999}

\title{Non-positively curved aspects of Artin groups\\of finite type}
\author{Mladen Bestvina}
\address{Department of Mathematics, University of Utah\\Salt Lake
City, UT 84112, USA}

\email{bestvina@math.utah.edu}

\newtheorem{thm}{Theorem}[section]
\newtheorem{lemma}[thm]{Lemma}
\newtheorem{cor}[thm]{Corollary}
\newtheorem{prop}[thm]{Proposition}
\newtheorem*{main}{Main Theorem}{}

\theoremstyle{remark}

\newtheorem{definition}[thm]{Definition}
\newtheorem{remark}[thm]{Remark}
\newtheorem{examples}[thm]{Examples}

\def\R{{\mathbb R}}
\def\Z{{\mathbb Z}}

\def\A{{\cal A}}
\def\G{{\cal G}}
\def\W{{\cal W}}
\def\P{{\cal P}}
\def\p{{\text {prod}}}

\def\Q{{\mathbb Q}}
\def\<{\langle}
\def\>{\rangle}

\begin{document}

\begin{abstract} Artin groups of finite type are not as well understood as
  braid groups. This is due to the additional geometric properties of braid
  groups coming from their close connection to mapping class groups. For each
  Artin group of finite type, we construct a space (simplicial complex)
  analogous to Teichm\" uller space that satisfies a weak nonpositive
  curvature condition and also a space ``at infinity'' analogous to the space of
  projective measured laminations. Using these constructs, we deduce several
  group-theoretic properties of Artin groups of finite type that are
  well-known in the case of braid groups.\end{abstract}




\maketitlepage

\section{Introduction and review} 
A theme in geometric group theory has been to, on one hand, construct
spaces with rich geometry where (given) interesting groups act by isometries,
and on the other hand, show that such actions have useful group-theoretic
consequences. A particularly successful sort of a geometric structure one
considers is a CAT(0) metric.  Experience shows that the class of groups that
act cocompactly and isometrically on CAT(0) spaces is large and includes many
standard groups, while at the same time the groups in this class have a rich
structure (see \cite{bh:cat0}). For example, the following properties hold in
this class:

\begin{itemize}
\item There are only finitely many conjugacy classes of finite
subgroups.
\item Every solvable subgroup is finitely generated and virtually
abelian.
\item The set of translation lengths of elements of infinite order is
bounded away from 0.
\end{itemize}

In this paper, we focus on the class of Artin groups of finite type
(and their quotients by the center). We review the basic definitions
and results in the next two sections. Tom Brady \cite{t.brady:artin}
has shown that Artin groups of finite type with $\leq 3$ generators
belong to the CAT(0) class and he has proposed a piecewise Euclidean
complex with a metric for each Artin group of finite type. Daan Krammer has
checked that this complex is CAT(0) for the braid group on 5 strands.

It is a result of Ruth Charney \cite{charney:biautomatic} that Artin
groups of finite type are biautomatic. This fact has a number of
group-theoretic consequences; in particular, by a result of
Gersten--Short \cite{gs:biautomatic}, the translation length of each
element of infinite order is positive, and consequently every
nilpotent subgroup is virtually abelian.

In this paper, we will construct a contractible simplicial complex on
which a given Artin group of finite type acts, and show that there is
a natural ``metric'' on it, preserved by the group action, that
satisfies a property somewhat weaker than CAT(0). The word ``metric''
is in quotes because it is not symmetric. This structure is, however,
sufficient to deduce the same group-theoretic properties as in the
case of CAT(0) metrics.

The main tool in this paper is the {\it left greedy normal form} for the
elements of Artin groups of finite type. The basic properties of this normal
form were established in the inspiring paper of Deligne \cite{deligne:fartin},
and also in Brieskorn--Saito \cite{bs:fartin}, both of which build on the work
of Garside \cite{garside:braids}. Charney \cite{charney:biautomatic} showed
that this form gives rise to a biautomatic structure. Deligne also showed that
the quotient by the Coxeter group of associated hyperplane complement is an
Eilenberg--Mac Lane space for each Artin group of finite type. Since these
manifolds can easily be compactified, Artin groups of finite type have finite
Eilenberg--Mac Lane spaces.

Aside from establishing the above group-theoretic properties of Artin
groups of finite type and their central quotients, we will give a
simple proof of Squier's theorem \cite{squier:artin} that these
groups are duality groups. Moreover, we show that they are highly
connected at infinity. (One way to show that a group is a duality
group of dimension $d$ is to argue that the group is $d$--dimensional
and $(d-2)$--connected at infinity.) We also show that every normal
abelian subgroup is central. This answers a question of Jim Carlson
that provided the impetus to study the geometry of Artin groups of
finite type.

We summarize the geometric features of the complex in the following
(see the next section for definitions of $S$ and $\Delta$):

\begin{main} Let $\cal A$ be an Artin group of finite type and
$\G=\A/\Delta^2$ the quotient of $\A$ by the central element
$\Delta^2$. Then $\G$ acts simplicially on a simplicial
complex $X=X(\G)$ with the following properties:
\begin{itemize}
\item The action is cocompact and transitive on the vertices {\rm(section
\ref{sec:defX})}.
\item $X$ is a flag complex {\rm(section
\ref{sec:defX})}.
\item $X$ is contractible, $(card (S)-3)$--connected at infinity and
proper homotopy equivalent to a $(card (S)-1)$--complex {\rm(Theorem
\ref{contractible})}.
\item Each pair of vertices of $X$ is joined by a preferred edge-path,
called a {\it geodesic} {\rm(section \ref{combing})}. The collection of
geodesics is invariant under the following operations: translation by
group elements, subpaths, inverses {\rm(Proposition \ref{symmetric})}, and
concatenations with nontrivial overlap.
\item There is a function $d=d_{wd}\co X^{(0)}\times X^{(0)}\to
\{0,1,2,\cdots\}$ (a ``non-symmetric distance function'') satisfying
the triangle inequality and
$$d(v,w)=0 \iff v=w$$ {\rm(Proposition \ref{symmetric})}. Further, $d_{wd}$
is comparable to the edge-path metric in $X$ {\rm(Lemma \ref{compare})}.
\item If $v_0,v_1,v_2,\cdots v_n$ are the consecutive vertices along a
geodesic, then $$d(v_0,v_n)=\sum_{i=1}^n d(v_{i-1},v_i)$$
\item If the associated Coxeter group $\W$ is irreducible and nonabelian, then
  every geodesic can be extended in both directions to a longer geodesic
  {\rm(Proposition \ref{charney})}.
\item {\rm(Theorem \ref{npc-D})} For any three vertices $a,b,c\in X$ and
any vertex $p\neq b,c$ on the geodesic from $b$ to $c$ we have
$$d(a,p)<\max\{d(a,b),d(a,c)\}.$$
\item The set of vertices of each simplex in $X$ has a canonical
cyclic ordering which is invariant under the group action and is
compatible with the passage to a face {\rm(Lemma \ref{ordering})}.
\end{itemize}
\end{main}

The dimension of the complex we construct is higher than the expected
dimension (ie, the dimension of the group). It is, however, often the
case that one has to increase the dimension in order to obtain better
geometric properties, cf non-uniform lattices, Teichm\" uller space,
Outer Space. 

Perhaps the most interesting remaining unresolved group-theoretic
questions about Artin groups $\A$ of finite type are the following:

\noindent{\bf Question 1}\qua Does $\A$ satisfy the Tits Alternative,
ie, if $H<\A$ is a subgroup which is not virtually abelian, does $H$
necessarily contain a nonabelian free group?

\noindent{\bf Question 2}\qua Is $\A$ virtually poly-free?

The Tits Alternative is not known for CAT(0) groups and it seems
unlikely that the techniques of this paper will resolve Question 1. A
group $G$ is {\it poly-free} if there is a finite sequence
$$G=G_0\supset G_1\supset G_2\supset\cdots\supset G_n=\{1\}$$ such
that each $G_i$ is normal in $G_{i-1}$ and the quotient $G_{i-1}/G_i$
is free. It would seem reasonable to expect that the pure Artin group
$P\A$ (the kernel of the homomorphism to the associated Coxeter group)
is poly-free, with the length of the series above equal to the
dimension (ie, the number of generators), and all successive
quotients of finite rank. This was verified by Brieskorn
\cite{brieskorn:bourbaki} for types $A_n$ (this is classical), $B_n$,
$D_n$, $I_2(p)$, and $F_4$. Of course, the positive answer to Question
2 implies the positive answer to Question 1. Another approach to
Question 1 would be to find a faithful linear representation of each
Artin group of finite type. This seems rather difficult (if not
impossible) even for braid groups.

Questions 1 and 2 can be asked in the setting of Artin groups of
infinite type as well. However, more basic questions are still open in
that setting; the most striking is whether all (finitely generated)
Artin groups admit a finite $K(\pi,1)$. Recently, Charney and Davis
\cite{cd:jams},\cite{cd:browder} have made substantial progress toward
this question, but the general case remains open.

We end the paper by constructing a space at infinity analogous to the
space of projectivized measured geodesic laminations in the case of
mapping class groups.

\rk{Acknowledgements}
I benefited from discussions about Artin groups with several people. I would
like to thank John Luecke and Alan Reid for having the patience to listen to
me while I was learning about Artin groups. Discussions with Warren Dicks have
shed much light on the difficulties in Question 2. I would also like to thank
the referee for suggesting simplifications of some of the arguments.

Support by the National Science Foundation is gratefully
acknowledged.

\subsection{Coxeter groups}

We follow the notation of \cite{cd:browder}. For proofs of the facts
listed below see \cite{bourbaki:coxeter}.  Let $S$ be a finite set. A
{\it Coxeter matrix} is a symmetric function $m\co S\times
S\to\{1,2,3,\cdots,\infty\}$ such that $m(s,s)=1$ for all $s\in S$ and
$m(s,s')\geq 2$ for $s\neq s'$. The associated {\it Coxeter group} is
the group $\W$ given by the presentation
$$\W=\<S\mid (ss')^{m(s,s')}=1\>$$ where $m=\infty$ means no relation.
If $S'$ is a proper subset of $S$ and $m'$ is the restriction of $m$ to
$S'\times S'$, then there is a natural homomorphism 
$$\W'=\<S'\mid (ss')^{m'(s,s')}=1\>\longrightarrow \W$$ and this homomorphism is
injective. The images of such homomorphisms are
called {\it special (Coxeter) subgroups} of $\W$. 

In this paper we will consider only {\it finite} Coxeter groups $\W$, and
this is what we assume from now on.

There is a
canonical faithful orthogonal representation of $\W$ on a finite-dimen\-sion\-al
vector space $V=V_\W$ (of dimension $\dim V=card(S)$) such that:
\begin{itemize}
\item Each generator $s\in S$ acts as a reflection, ie, it fixes a
  codimension 1 subspace of $V$. These subspaces and their
  $\W$--translates are called {\it walls}.
\item The closures in $V$ of the complementary components of the union
  of the walls are the
{\it chambers}. They are simplicial cones, and $\W$ acts on the set of
  chambers simply transitively.
\item Special subgroups and their conjugates are precisely the
  stabilizers of nonzero points in $V$.
\item There is a chamber $Q$, called the {\it fundamental chamber},
  that intersects the fixed sets of each $s\in S$ in top-dimensional
  faces. When $\W$ is irreducible, $Q$ and $-Q$ are
  the only such chambers.
\item The {\it longest element} is the unique element $\Delta\in\W$
  that takes $Q$ to $-Q$ (and {\it vice versa}). It has order 2, and
  conjugation by $\Delta$ induces an involution of $S$.
\end{itemize}

Two chambers $Q'$ and $Q''$ are {\it adjacent} if $Q'\cap Q''$ is
contained in a wall and has nonempty interior in the wall. This wall
is then the unique wall that separates $Q'$ from $Q''$, and we say
that it {\it abuts} $Q'$ and $Q''$. A sequence $Q_0,Q_1,\cdots, Q_n$
of chambers is a {\it gallery $G$ of length $n$ (from $Q_0$ to $Q_n$)}
if $Q_i$ and $Q_{i+1}$ are adjacent for $i=0,1,2,\cdots,n-1$. Any path
in $U$ which is transverse to the collection of the walls determines a
gallery and every gallery arises in this way. Any two chambers are
connected by a gallery. The {\it distance} between two
chambers $A,B$ is the length of a shortest gallery connecting them;
equivalently, it is the number of walls separating $A$ and $B$. A
gallery is {\it geodesic} if its length is equal to the distance
between the initial and the terminal chamber. Equivalently, a gallery
is geodesic if it crosses each wall at most once (we say that a
gallery crosses a wall if any associated path does, ie, if the wall
separates a pair of consecutive chambers). For example, straight line
segments transverse to the collection of walls determine geodesic
galleries. 

A Coxeter group $\W$ is {\it irreducible} if there is no nontrivial
partition $S=S_1\sqcup S_2$ with $m(s_1,s_2)=2$ whenever $s_1\in S_1$
and $s_2\in S_2$.

\subsection{Artin groups}

If $S$ is a finite set and $m\co S\times S\to \{2,3,\cdots,\infty\}$ a
Coxeter matrix, the associated {\it Artin group} is defined to be the
group
$$\A=\<S\mid \p(s,s',m(s,s'))=\p(s',s,m(s',s))\>$$
where $$\p(s,s',m)=\underbrace{ss'ss'\cdots}_{m\text{\ times}}$$ and
again $m=\infty$ means no relation.

There is a natural homomorphism $\pi\co \A\to\W$, $\pi(s)=s$. The kernel
is the {\it pure Artin group} $\P$. If $S'$ is a subset of $S$, let
$\A'$ be the Artin group associated with the restriction of $m$ to
$S'\times S'$. We then have a natural homomorphism $\A'\to\A$,
$s'\mapsto s'$. This homomorphism is injective.

By $\A^+$ we denote the submonoid of $\A$ generated by $S$, ie, the
monoid of positive words. As shown by Deligne and Brieskorn--Saito, if
two positive words represent the same element of $\A$, then they can
be transformed to each other by repeated substitutions given by the
defining relations. A positive word $s_{i_1}s_{i_2}\cdots s_{i_k}$ can
be geometrically thought of as a gallery of length $k$ starting at
$Q$: $Q,s_{i_1}(Q),s_{i_1}s_{i_2}(Q),\cdots,s_{i_1}s_{i_2}\cdots
s_{i_k}(Q)$. We will also consider any $\W$--translate of this gallery
as being associated with the same positive word. Multiplication in
$\A^+$ corresponds to the concatenation of galleries (after a possible
$\W$--translation so that the last chamber of the first gallery
coincides with the first chamber of the second gallery). 

\subsection{Normal form in $\A^+$ and $\A$}

Following Deligne, for $x,y\in\A^+$ we write $x<y$ if there is
$z\in\A^+$ such that $y=xz$. We say that $A\in\A^+$ is an {\it atom}
if its word-length in $\A^+$ is equal to the word-length of
$\pi(A)\in\W$, or equivalently if a gallery associated to $A$ is geodesic.
The homomorphism $\pi\co \A\to\W$ induces a
bijection between the set of atoms and $\W$. We will often identify an
atom $A$ with the associated chamber $\pi(A)(Q)$. The atom that corresponds
to the longest element $\Delta$ is also denoted $\Delta$. We have
$1<A<\Delta$ for all atoms $A$. $\Delta^2$ is central in $\A$ and
conjugation by $\Delta$ induces an involution of $\A$ that extends the
involution $S\to S$ mentioned above. This involution will be denoted
by $x\mapsto\overline x$. It restricts to involutions of $\A^+$ and of
the set of atoms.

The following two propositions are proved in Deligne's paper
\cite{deligne:fartin}.

\begin{prop} \label{characterization}
Suppose that $\Omega$ is a nonempty finite subset of $\A^+$ such
that \begin{enumerate} \item $x\in\Omega$, $x'<x$ implies
$x'\in\Omega$, \item if $\sigma,\tau\in S$ are two generators, $x\in\A^+$,
$x\sigma,x\tau\in\Omega$, then
$$x\text{\p}(\sigma,\tau,m(\sigma,\tau))\in\Omega.$$ \end{enumerate} 

Then there is a unique $y\in\A^+$ such that $$\Omega=
\{ x\in\A^+\mid x<y\}.$$ 

Further, for any $y$ the set $\Omega=\{ x\in\A^+\mid x<y\}$ satisfies 1 and 2.
\qed\end{prop}

In particular, for $x_1,x_2\in\A^+$ we can define $x_1\wedge
x_2\in\A^+$ as the largest element $z\in\A^+$ such that $z<x_1$ and
$z<x_2$. The existence of the largest such element follows by applying
Proposition \ref{characterization} to the set
$\Omega=\{z\mid z<x_1,z<x_2\}$. 

There is a function $reverse\co \A^+\to\A^+$ given by
$reverse(s_{i_1}s_{i_2}\cdots s_{i_k})=s_{i_k}\cdots
s_{i_2}s_{i_1}$. If $x\in\A^+$ is represented by a gallery
$Q_0,Q_1,\cdots,Q_{k-1},Q_k$, then $reverse(x)$ is represented by the gallery
$Q_k,Q_{k-1},\cdots,Q_1,Q_0$. 

For every atom $A$ there is an atom $A^*$ such that
$AA^*=\Delta$ (and there is also an atom $^*A$ such that $^*AA=\Delta$). 
Every $x\in\A^+$ can be represented as the product of atoms (or even
generators), say $x=A_1A_2\cdots A_k$. Then we have 
$$xA_k^*\overline A_{k-1}^* A_{k-2}^*\cdots=\Delta^k$$ so that
$x<\Delta^k$. Further, if $x,y<\Delta^k$, and if $\tilde x,\tilde
y\in\A^+$ are such that $x\tilde x=\Delta^k$ and $y\tilde y=\Delta^k$,
then $x<y\iff reverse(\tilde y)<reverse(\tilde x)$ and we deduce (from
the existence of $\tilde x\wedge \tilde y$) that for any $x,y\in\A^+$
there is a unique $x\vee y\in\A^+$ such that \begin{itemize}
\item $x<x\vee y$, \item $y<x\vee y$, and \item $x<z$, $y<z$
implies $x\vee y<z$.\end{itemize}

For every $x\in\A$ there is $k\geq 0$ such that $\Delta^k
x\in\A^+$. This reduces the understanding of $\A$ to the understanding
of $\A^+$. There is no such reduction when the associated Coxeter
group is infinite.

\begin{prop} \label{nf} For every $g\in\A^+$ there is a unique atom $\alpha(g)$
  such that \begin{enumerate} \item $\alpha(g)<g$, and \item if $A$ is
  an atom with $A<g$, then $A<\alpha(g)$.\end{enumerate} Furthermore,
  $\alpha(xy)=\alpha(x\alpha(y))$ for all $x,y\in\A^+$.\qed
\end{prop}

The {\it left greedy normal form}
of $g\in\A^+$ is the representation of $g$ as the product of
nontrivial atoms
$$g=A_1A_2\cdots A_k$$
such that $A_i=\alpha(A_iA_{i+1}\cdots A_k)$
for all $1\leq i\leq k$. To emphasize that the above product of atoms
is the normal form of $g$, we will write $g=A_1\cdot
A_2\cdot\cdots\cdot A_k$. Clearly, the normal form for $g$ is unique.
A product $A_1A_2\cdots A_k$ is a normal form iff each subproduct of
the form $A_iA_{i+1}$ is a normal form. If $A_1\cdot
A_2\cdot\cdots\cdot A_k$ is a normal form, then so is $\overline
A_1\cdot \overline A_2\cdot\cdots\cdot \overline A_k$.

$\Delta$'s appear at the beginning of the normal form, ie,
if $A_1\cdot A_2\cdot\cdots\cdot A_k$ is a normal form, then there is
$j$ such that $A_1=A_2=\cdots=A_j=\Delta$ and $A_i\neq\Delta$ for
$i>j$. 

Geometric interpretation of the normal form is given by the following
result.

\begin{prop} \label{interpretation}
Let $A$ and $B$ be nontrivial atoms, with $A$ represented by a
  geodesic gallery from $Q_1$ to $Q_2$ and $B$ by a geodesic gallery
  from $Q_2$ to $Q_3$. Then $A\cdot B$ is a normal form
  iff no wall that abuts $Q_2$ has $Q_1$ and $Q_2$ on one side and
  $Q_3$ on the other side.\qed
\end{prop} 

There is the homomorphism $length\co \A\to\Z$ that sends each $s\in S$ to
$1\in\Z$. This homomorphism restricts to the word-length on $\A^+$.

We will need a straightforward generalization of Proposition \ref{nf}.

\begin{prop} \label{summary+}
Suppose $x\in\A^+$ and $y_1\cdot y_2\cdot\cdots\cdot y_l$ is a normal
form in $\A^+$. Suppose the product 
$x(y_1\cdot y_2\cdot \cdots\cdot
y_l)$
has normal form $$z_1\cdot z_2\cdot\cdots\cdot z_m.$$
Then for each $i\leq l$ the normal form of the product
$x(y_1y_2\cdots y_i)$ begins with $$z_1\cdot
z_2\cdot\cdots \cdot z_i\cdot\cdots$$
\end{prop}
 
\begin{proof}
For $i=1$ this follows from Proposition \ref{nf}.
We argue by induction on $i$. It suffices to show that $z_1z_2\cdots
z_i<xy_1y_2\cdots y_i$. By Proposition \ref{nf} we have
$z_1<xy_1$ so we can write $xy_1=z_1w$ for some $w\in\A^+$. Then
$xy_1y_2\cdots y_l=z_1z_2\cdots z_m$, after cancelling on the left,
implies that $wy_2\cdots y_l=z_2\cdots z_m$ and the right-hand side is
visibly a normal form. By induction it follows that $z_2\cdots
z_i<wy_2\cdots y_i$. Now multiply on the left by $z_1$ to obtain
$z_1z_2\cdots z_i<xy_1y_2\cdots y_i$.
\end{proof}

\section{$\G$ and its complex $X(\G)$}

\subsection{Definition of $X(\G)$}\label{sec:defX}

It is more convenient to study the group $\G=\A/\Delta^2$. The groups
$\G\times\Z$ and $\A$ are commensurable. Indeed, the
homomorphism $\A\to\G\times\Z$ which is natural projection in the
first coordinate and sends each generator to $1\in\Z$ in the second
coordinate (ie, it is the length homomorphism in the second
coordinate) is a monomorphism onto a finite index subgroup. Moreover,
$\Z=\<\Delta^2\>\hookrightarrow \A$ has a splitting with values in
$\frac{1}{length(\Delta^2)}\Z$ (given by length divided by the length
of $\Delta^2$).

The goal of this section is to describe a contractible simplicial
complex $X=X(\G)$ on which the group $\G=\A/\<\Delta^2\>$ acts cocompactly and
with finite point stabilizers. The vertex set is the coset space
$$V=\A/\<\Delta\>=\{ g\<\Delta\>\mid g\in \A\}.$$ The group $\A$ acts naturally
on the left, and the central element $\Delta^2$ acts trivially, so
$\G=\A/\<\Delta^2\>$ acts on $V$. The action is transitive, and the
stabilizer of the point $\<\Delta\>\in V$ is $\<\Delta\>\cong \Z/2$.

We call the elements of $V$ the {\it vertices}. Each vertex has a
representative $g\in \A^+$ (obtained by multiplying an arbitrary
representative by a high power of $\Delta^2$). Let $g=\Delta^m\cdot
B_1\cdot B_2\cdot \cdots \cdot B_k$ be the normal form for $g$ with
$B_i\neq\Delta$ for all $i$ (as remarked above,
if the normal form of $g$ involves any $\Delta$'s, they appear at the
beginning). Now shift all $\Delta$'s to the end to obtain a
representative in $\A^+$ of the same coset whose normal form does not
have any $\Delta$'s. If $m$ is even, this amounts to replacing $g$ by
$B_1\cdot B_2\cdot \cdots \cdot B_k$, and if $m$ is odd $g$ is
replaced by $\overline B_1\cdot \overline B_2\cdot \cdots \cdot
\overline B_k$. Such a coset representative is unique, and we will
frequently identify $V$ with the set of elements of $\A^+$ whose
normal form has no $\Delta$'s. The identity element of $\A^+$ (whose
normal form is empty) is viewed as the basepoint, denoted $*$.

The {\it atomnorm} of a vertex $v$, denoted $|v|$, is the number of
atoms in the normal form of the special representative, ie, the
number of atoms not counting $\Delta$'s in the normal form of any
positive representative.  Left translation by $\Delta$ sends a special
representative $B_1\cdot B_2\cdot\cdots \cdot B_m$ to a special
representative $\overline B_1\cdot \overline B_2\cdot\cdots \cdot
\overline B_m$ so this involution is atomnorm preserving. It follows that
there is a unique left invariant function $d_{at}\co V\times
V\to\Z_+=\{0,1,2,\cdots\}$ (ie $d_{at}(g(v),g(w))=d_{at}(v,w)$ for all $g\in
\G$ and all $v,w\in V$) such that $d_{at}(*,v)$ is the atomnorm of $v$.
We call $d_{at}$ the {\it atomdistance} on $V$.

In order to prove that $d_{at}$ is symmetric and satisfies the triangle
inequality, we need the following lemma.

\begin{lemma} \label{l1}
If $g\in \A^+$ is the product of $k$ atoms, then the
normal form of $g$ has $\leq k$ atoms.\end{lemma}

\begin{proof} Induction on $k$ starting with $k=1$ when it's clear. 
For $k=2$ the statement follows from the definition of normal form and
the fact that a subgallery of a geodesic gallery is a geodesic
gallery. Say $g=B_1B_2\cdots B_k$. Inductively, the normal form for
$B_2\cdots B_k$ is $C_2\cdot C_3\cdot\cdots\cdot C_l$ with $l\leq
k$. If $B_1C_2$ is an atom, then the normal form for $g$ is
$(B_1C_2)\cdot C_3\cdot\cdots\cdot C_l$ and has length $l-1<k$. If
$B_1C_2$ is not an atom, then the normal form for $B_1C_2$ is
$X\cdot Y$ (say) and the normal form for $g$ is $X$ followed by the
normal form for $YC_3\cdots C_l$. The latter has length $<l$ by
induction, so the statement is proved. \end{proof}

\begin{prop}
$d_{at}\co V\times V\to\Z_+$ is a distance function.
\end{prop}

\begin{proof} The triangle inequality follows from Lemma \ref{l1}. 
We argue that $d_{at}$ is symmetric. By left invariance, it suffices
to show that $d_{at}(*,v)=d_{at}(v,*)$ for all $v\in V$.  Let $g$ be a
group element such that $g(v)=*$; we need to argue that the norm of
$v$ equals the norm of $g(*)$ (which can be viewed as an ``inverse''
of $v$; there are two inverses --- obtained from each other by applying
the ``bar'' involution). Say $v=B_1\cdot B_2\cdot \cdots\cdot B_k$ is
the normal form (without $\Delta$'s). Let $C_i$ be the unique atom
such that $C_iB_i=\Delta$. Thus $\Delta C_i\in\G$ is the inverse of
$B_i$ and $g$ can be taken to be $\Delta C_k\Delta C_{k-1}\cdots
\Delta C_1$. Thus $g(*)$ is represented by $\overline C_k
C_{k-1}\overline C_{k-2}\cdots$, a product of $k$ atoms. By Lemma \ref{l1},
the norm of $g(*)$ is $\leq k$. Applying the same argument with the
roles of $v$ and $g(v)$ reversed, we see that the norms of $v$ and
$g(v)$ are equal. \end{proof}

\begin{definition} $X(\G)$ is the simplicial complex whose vertex set
is $V$ and a collection $v_0,v_1,\cdots,v_m\in V$ spans a simplex iff
$d_{at}(v_i,v_j)=1$ for all $i\neq j$.\end{definition}

Thus $X(\G)$ is a flag complex, ie, if each pair in a set of vertices
bounds an edge in $X(\G)$ then this set spans a simplex.
The group $\G$ acts on $X(\G)$ simplicially and cocompactly. Note that
the distance $d_{at}$ on $V$ is the edge-path distance between the
vertices of $X(\G)$.

\begin{examples} We can explicitly see $X(\G)$ in simple cases. When
$\A=\Z\times \Z$, $X(\G)$ is the line. When $\A$ is $\Z\times \Z\times
\Z$, $X$ is the plane (triangulated in the usual $(3,3,3)$--fashion),
and more generally, when $\A=\Z^m$, then $X(\G)$ is Euclidean
$(m-1)$--space. When $\A=B_3$ is the braid group on 3 strands, $X(\G)$ is
the union of triangles glued to each other along their vertices so
that the spine is the trivalent tree illustrated in figure 1.\end{examples}

\begin{figure}[ht!]
\centerline{
\relabelbox\footnotesize
\epsfxsize 0.9\textwidth\epsfbox{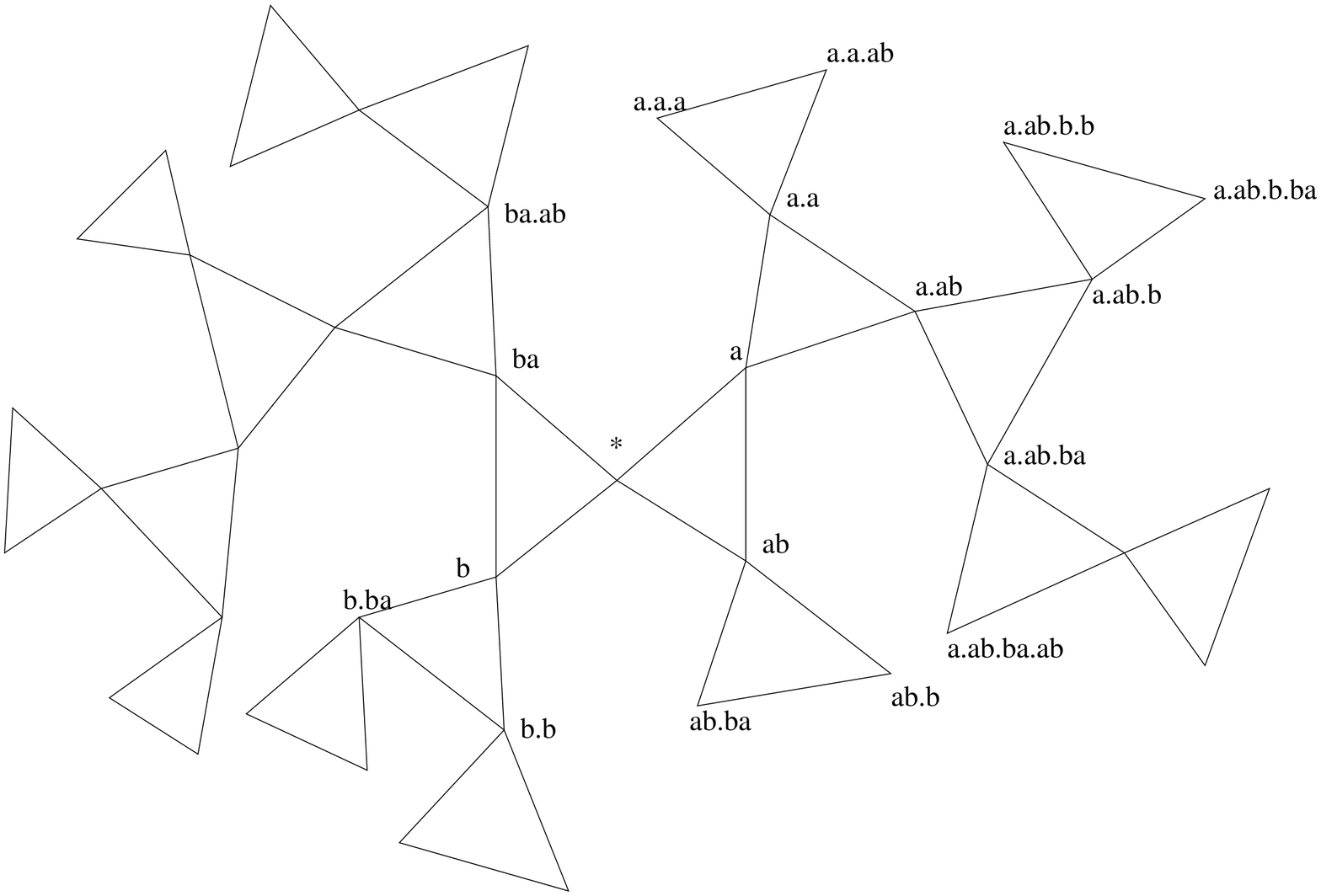}
\relabela <-2pt,0pt> {a}{$a$}
\relabel {a.a}{$a.a$}
\relabela <-2pt,2pt> {a.a.a}{$a.a.a$}
\relabel {a.a.ab}{$a.a.ab$}
\relabel {a.ab}{$a.ab$}
\relabel {a.ab.b}{$a.ab.b$}
\relabel {a.ab.ba}{$a.ab.ba$}
\relabel {a.ab.ba.ab}{$a.ab.ba.ab$}
\relabel {a.ab.b.b}{$a.ab.b.b$}
\relabel {a.ab.b.ba}{$a.ab.b.ba$}
\relabel {ab}{$ab$}
\relabel {ab.b}{$ab.b$}
\relabel {ab.ba}{$ab.ba$}
\relabela <0pt,-2pt> {*}{$*$}
\relabel {b}{$b$}
\relabel {ba}{$ba$}
\relabel {b.b}{$b.b$}
\relabel {ba.ab}{$ba.ab$}
\relabela <-2pt,2pt> {b.ba}{$b.ba$}
\endrelabelbox
}
\caption{$X(\G)$ for $\G=\A/\Delta^2$, $\A=\<a,b\mid aba=bab\>$}
\end{figure}

\subsection{Related complexes} \label{sec:related}
One can build a similar complex $X'(\G)$
with the vertex set $\G=\A/\Delta^2$. Each vertex has a special
representative whose normal form either has no $\Delta$'s, or has a
single $\Delta$, and in the latter case we agree to push this $\Delta$
to the last slot. There is an oriented edge from a vertex $v$ to a
vertex $w$ if the special representative of $w$ is obtained from the
special representative of $v$ by rightmultiplying by an atom. The
simplicial complex $X'(\G)$ is defined to be the flag complex
determined by the resulting graph, ie, a collection of vertices spans
a simplex provided all pairs span an edge. The natural quotient map
$\G\to\G/\<\Delta\>$ extends to a simplicial map $p\co X'(\G)\to X(\G)$. The
preimage of the basepoint $*\in X(\G)$ is the edge $[*,\Delta]\subset
X'(\G)$. More generally, the preimage of a simplex $\sigma$ is the
(triangulated) prism $\sigma\times [0,1]$. There is a global homeomorphism
$X'(\G)\cong X(\G)\times [0,1]$. In the first coordinate this
homeomorphism is given by $p$ and in the second it is the simplicial
map that sends the vertices whose special representative has no
$\Delta$'s to $0$ and the vertices whose special representative has
one $\Delta$ to 1.

Analogously, there is a complex $X(\A)$ with vertex set $\A$
and similarly defined simplicial structure:
edges are drawn from $*$ to the atoms and extended equivariantly,
and then the higher-dimensional simplices are filled in. The quotient
map extends to the natural simplicial map
$q\co X(\A)\to X(\G)$ and there is a global homeomorphism $X(\A)\cong
X(\G)\times \R$ in the first coordinate given by $q$ and in the second
it is the map that is linear on each simplex and sends a vertex $g\in
\A$ to $length(g)\in\Z$ (recall that $length\co \A\to\Z$ is the
homomorphism that sends each $s\in S$ to $1\in\Z$).
The fact that $\Delta^2$
virtually splits is transparent in this model.

We end this discussion by recalling (see \cite{cd:browder}) the description of
another important $\A$--complex and sketching the proof of its
contractibility using the contractibility of $X(\A)$ (Theorem
\ref{contractible}). 

The Coxeter sphere $\Sigma_{\W}$ can be described as (the geometric
realization of) the poset of cosets $w\W_F$ where $w\in\W$ and $\W_F$ is the
special Coxeter subgroup of $\W$ generated by $F\subset S$, as $F$ runs over
all proper subsets (including $\emptyset$) of $S$. The partial order is given
by inclusion. Here $\W_F$ corresponds to the barycenter of the largest face of
$Q$ fixed by $F$ and $w\W_F$ corresponds to the $w$--translate of this face.
One description \cite{cd:browder} of the {\it universal cover of the Salvetti
  complex} $\tilde \Sigma_{\W}$ associated to $\A$ (and $\W$) is that it is
the poset of subsets $a\W_F\subset\A$ with $a\in\A$ and $F$ any subset of $S$
(including $\emptyset$ and $S$), where $\W\supset\W_F$ is identified with the
set of all atoms in $\A$.  The Salvetti complex is the quotient of $\tilde
\Sigma_{\W}$ by the pure Artin group. The Salvetti complex is homotopy
equivalent to the associated hyperplane complement (see \cite{cd:browder}).

The cover of $\tilde \Sigma_{\W}$ by the translates of $\W_S$ has the property
that all nonempty intersections are contractible (this follows
from the existence of least upper bounds for sets of atoms). Therefore,
$\tilde \Sigma_{\W}$ is homotopy equivalent to the nerve of this cover.

Now consider the cover of $X(\A)$ by the translates of the subcomplex spanned
by all atoms. This cover also has all nonempty intersections contractible and
its nerve is isomorphic to the nerve above. Thus $\tilde \Sigma_{\W}$ and
$X(\A)$ are homotopy equivalent.

\section{Geometric properties of $X(\G)$}

\subsection{Geodesics}\label{combing}

There is a canonical ``combing'' of $X$: for each vertex $v\in X$ we
have a canonical edge-path from the basepoint $*$ to $v$. If
$v=B_1\cdot B_2\cdot\cdots\cdot B_k$ is the normal form without
$\Delta$'s, then the edge-path associated to $v$ is $$*,B_1,B_1\cdot
B_2,\cdots,B_1B_2\cdots B_k=v.$$ Note that left multiplication by
$\Delta$ sends $v$ to $\overline B_1\cdot \overline B_2\cdots
\overline B_k$ and it sends the combing path from $*$ to $v$ to the
combing path from $*=\Delta(*)$ to $\Delta(v)$. It follows that by
left translating we obtain a canonical path from any vertex to any
other vertex. We call all such edge-paths {\it geodesics}.

It is interesting that the geodesics are symmetric. This fact
follows from the work of R Charney (Lemma 2.3 of
\cite{charney:symmetric}). For
completeness we indicate a proof.

\begin{prop}\label{symmetric}
  The geodesic from $v$ to $w$ is the inverse of the geodesic
  from $w$ to $v$.
\end{prop}

\begin{proof} First note that an edge-path is a geodesic if and
  only if every subpath of length 2 is a geodesic. Thus it
  suffices to prove the Proposition in the case that the geodesic
  from $v$ to $w$ has length 2. Further, after left-translating, we
  may assume that the middle vertex of the path is $*$, and so $v$ and
  $w$ are atoms, say $A$ and $B$ respectively. Denote by $C$ the atom
  with $CA=\Delta$. Note that $C$ is represented by a geodesic gallery
  from $-A$ to $Q$. We now have (recall our convention that we
  identify an atom $A$ with the chamber $\pi(A)(Q)$):
  
  $A,*,B$ is a geodesic $\iff$ $CA=\Delta=*,C,CB$ is a geodesic
  $\iff$ $C\cdot B$ is a normal form $\iff$ (by Proposition
  \ref{interpretation}, see figure 2) no wall $W=Fix(s)$ that abuts
  $Q$ has $Q$ and $-A$ on one side and $B$ on the other $\iff$
  (applying ``$-$'') no wall that abuts $-Q$ has $-Q$ and $A$ on one
  side and $-B$ on the other $\iff$ (using that a wall abuts $Q$ iff
  it abuts $-Q$ and it always separates the two) no wall that abuts
  $Q$ has $Q$ and $-B$ on one side and $A$ on the other.
  
\begin{figure}[ht!]
\centerline{
\relabelbox\small
\epsfxsize0.5\textwidth\epsfbox{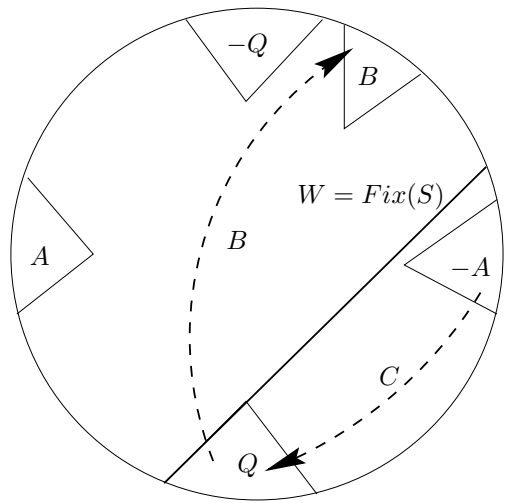}
\relabel {A}{$A$}
\relabel {-A}{$-A$}
\relabel {B}{$B$}
\relabel {B1}{$B$}
\relabel {C}{$C$}
\relabel {Q}{$Q$}
\relabel {-Q}{$-Q$}
\relabela <-2pt,15pt> {W}{$W=Fix(S)$}
\endrelabelbox}
\caption{If $W$ exists, $C\cdot B$ is not a normal form; $C$ can be
  extended to a geodesic gallery $C'=Cs$ with $C'<CB$.}
\end{figure}

  But the last statement is exactly the same as a previous statement
  with the roles of $A$ and $B$ interchanged.
\end{proof}

Recall that if $A$ is an atom, the atom $A^*$ is defined by $AA^*=\Delta$. The
following fact is contained in \cite{charney:symmetric} (cf Lemmas 2.3 and
2.5). 

\begin{cor}\label{new}
If $A_1\cdot A_2\cdot \cdots\cdot A_k$ is a normal form, so are $A_k^*\cdot
\overline A_{k-1}^*\cdot A_{k-2}^*\cdot\cdots$.
\end{cor}

\begin{proof}
  Left-translate the path $Q,A_k^*,A_k^*\overline A_{k-1}^*,\cdots$ by
  $A_1A_2\cdots A_k$ to get the path $A_1A_2\cdots A_k,A_1A_2\cdots
  A_{k-1},\cdots,A_1,Q$. Since the inverse of this path is a geodesic by
  assumption, so is the original path by Proposition
  \ref{symmetric}.
\end{proof}

\subsection{Worddistance}\label{sec:worddistance} 

By the {\it wordnorm} of a vertex $v\in V$, denoted
$||v||$, we mean the word-length of the special representative of $v$,
or equivalently, the word-length of the normal form of any
representative of $v$ after discarding the $\Delta$'s. Observe that
left translation by $\Delta$ is wordnorm preserving (it sends $B_1\cdot
B_2\cdot\cdots\cdot B_m$ to $\overline B_1\cdot
\overline B_2\cdot\cdots\cdot \overline B_m$), and so there is a
unique function $d_{wd}\co V\times V\to \Z_+$ that is left-invariant and
$d_{wd}(*,v)=||v||$ for all $v\in V$. We call $d_{wd}$ the {\it worddistance},
but we caution the reader that $d_{wd}$ is not symmetric. 
For example, if $B$ is an atom, then $d_{wd}(*,B)$ is
the word-length of $B$, while $d_{wd}(B,*)$ is the word-length of the
atom $C$ complementary to $B$ (ie such that $CB=\Delta$). This
observation generalizes as follows.

\begin{lemma} \label{compare}
  $d_{wd}(v,w)+d_{wd}(w,v)=d_{at}(v,w)\delta$, where $\delta$ is the
  word-length of $\Delta$.
\end{lemma}

\begin{proof} After left-translating, we may assume that $v=*$. If
  $w=B_1\cdot B_2\cdot\cdots\cdot B_k$ (without $\Delta$'s), then
  $d_{at}(*,w)=k$ and $d_{wd}(*,w)$ is the word-length of $B_1\cdots B_k$. To
  compute $d_{wd}(w,*)$, denote by $C_i$ the atom with $C_iB_i=\Delta$. Left
  translate by $\cdots C_{3}\overline C_{2}C_1$ to get
  $d_{wd}(w,*)=d_{wd}(\Delta^k,\cdots C_{3}\overline C_{2}C_1)=||\cdots
  C_{3}\overline C_{2}C_1||$. Since the word-length of $C_i$ is obtained from
  $\delta$ by subtracting the word-length of $B_i$, the word-length of $\cdots
  C_{3}\overline C_{2}C_1$ equals $k\delta-d_{wd}(*,w)$, so it remains to
  argue that $\cdots C_{3}\cdot \overline C_{2}\cdot C_1$ is a normal
  form. But this is an immediate consequence of Corollary \ref{new}
  ($C_i=\overline B_i^*$).
\end{proof}

The triangle inequality $d_{wd}(v,w)+d_{wd}(w,z)\geq d_{wd}(v,z)$ is
immediate, and so are the inequalities $d_{at}(v,w)\leq
d_{wd}(v,w)\leq \delta d_{at}(v,w)$, $d_{wd}(w,v)\leq \delta
d_{wd}(v,w)$ (showing that $d_{wd}$ is ``quasi-symmetric'', and that
$d_{wd}$ and $d_{at}$ are quasi-isometric).

We will next argue that the wordnorms of adjacent vertices are always
distinct. 
To this end, we examine the relationship between the normal forms
associated to two adjacent vertices, say $v=B_1\cdot
B_2\cdot\cdots\cdot B_k$ and $w=C_1\cdot C_2\cdot\cdots\cdot C_l$
(without $\Delta$'s). Since $d_{at}(v,w)=1$, there is an atom $X\neq
1,\Delta$ such that $vX$ and $w$ represent the same coset,
ie $$vX=w\Delta^m$$ for some integer $m$. Similarly, there is a
atom $Y\neq 1,\Delta$ such that
$$wY=v\Delta^n$$ for an integer $n$.

\begin{lemma} One of $m,n$ is 0 and the other is 1.\end{lemma}

\begin{proof} Substituting $v=wY\Delta^{-n}$ into $vX=w\Delta^m$ yields
  $Y\Delta^{-n}X=\Delta^m$.  If $n$ is even, this equation implies
  that $YX$ is a power of $\Delta$, hence $YX=\Delta$ and $m+n=1$. If
  $n$ is odd, then similarly $\overline YX=\Delta$ and again $m+n=1$.
  
  It remains to argue that $m,n\geq 0$. Assume $n<0$. Then
  $v=wY\Delta^{-n}$ will have $\Delta$'s in its normal form (as
  $\Delta$ can be moved to the front), contradiction. Similarly $m\geq
  0$.
\end{proof}

It now follows that either $||v||<||w||$ (if $m=0$) or $||w||<||v||$
(if $n=0$).  We orient the edge $[v,w]$ from $v$ to $w$ if
$||v||<||w||$. Thus the arrows point away from the basepoint. Note
that these orientations are not group invariant.  

It is clear from the definition that there are no oriented cycles. This
implies that the dimension of $X(\G)$ is $length(\Delta)-1$. Indeed, the set
of vertices of a simplex of $X(\G)$ is totally ordered and thus a
top-dimensional simplex that, say, contains the base point $*$ corresponds to
a chain of atoms $*<B_1<B_2<\cdots <B_m$ of maximal possible length
$length(\Delta)$.

The following lemma sharpens Lemma \ref{l1}.

\begin{lemma} \label{l4}
Suppose $g=X_1X_2\cdots X_k\in \A^+$ is written as the product of $k$
atoms, and let $g=B_1\cdot B_2\cdot\cdots\cdot
B_l$ be the normal form for $g$. Then $l\leq k$ and 
$X_1X_2\cdots X_i<B_1B_2\cdots
B_i$ for all $1\leq i\leq l$.\end{lemma}

\begin{proof} That $l\leq k$ is the content of Lemma \ref{l1}. 
That $X_1<B_1$ follows from the definition of the normal form. Assume
inductively that $X_1\cdots X_{i-1}<B_1\cdots B_{i-1}$. Thus we can
write $B_1\cdots B_{i-1}=X_1\cdots X_{i-1}y$ for some $y\in \A^+$, and thus
$X_iX_{i+1}\cdots X_k=yB_iB_{i+1}\cdots B_l$. But then
$$X_i<\alpha(yB_iB_{i+1}\cdots B_l)=\alpha(yB_i)<yB_i$$ so that
multiplication on the left by $X_1\cdots X_{i-1}$ yields $X_1\cdots
X_i<B_1\cdots B_i$. The equality in the displayed formula follows from
Proposition \ref{nf}. \end{proof}

\subsection{Contractibility of $X(\G)$ and its topology at infinity}
\label{sec:duality}

Denote by $d=card (S)-1$ the dimension of the Coxeter sphere (the unit
sphere in $V_\W$).

In this section we prove:

\begin{thm} \label{contractible} 
  $X(\G)$ is contractible and it is proper homotopy equivalent to a cell
  complex which is $d$--dimensional and $(d-2)$--connected at infinity.
\end{thm}

In particular, this recovers a theorem of C\,C Squier
  \cite{squier:artin}: $\G$ is a virtual duality group of dimension $d$
  and $\A$ is a duality group of dimension $d+1$.

The method of proof is to consider the function $F\co X\to [0,\infty)$
that to a vertex $v$ assigns the wordnorm $F(v)=||v||$ and is linear on each
simplex. The theorem follows by standard methods (see \cite{bb:morse},
\cite{bf:duality}) from:

\begin{prop} \label{links} $F$ is nonconstant on each edge of $X(\G)$. The
  ascending link at each vertex $v$ is either contractible or homotopy
  equivalent to the sphere $S^{d-1}$ and the descending link at each
  vertex $v\neq *$ is contractible.
\end{prop}

For a vertex $v$, the {\it ascending link} $Lk_\uparrow(v)$ [{\it
descending link} $Lk_\downarrow(v)$] is defined to be the link of $v$ in
the subcomplex of $X(\G)$ spanned by the vertices $w$ with $F(w)\geq
F(v)$ [$F(w)\leq F(v)$].

We need a lemma. The proof was suggested by the referee and it is considerably
simpler than the original argument.

\begin{lemma} \label{below}
Given an atom $A$ and an element $b\in \A^+$ there is
a unique atom $C$ such that for every $u\in \A^+$
$$A<bu \iff C<u.$$
\end{lemma}

\begin{proof} 
Let $d=A\vee b$. Then $d=bc$ for some $c\in \A^+$ and
$$A<bu\Leftrightarrow d<bu\Leftrightarrow c<u.$$ Taking $u=\Delta$, we see
that $c$ is an atom.
\end{proof}

\begin{proof}[Proof of Proposition \ref{links}]
  Represent $v$ as a normal form $B_1\cdot B_2\cdot\cdots\cdot B_k$ without
  $\Delta$'s. A representative of an adjacent vertex can be obtained from
  $B_1\cdots B_k$ by multiplying on the right by an atom, say $B$. Note that
  $F(vB)<F(v)$ precisely when $\Delta<B_1\cdots B_kB$. By Lemma \ref{below},
  there is an atom $C$ such that for any atom $B$ $$\Delta<B_1\cdots B_kB\iff
  C<B.$$
  Thus the descending link at $v$ can be identified with the poset of
  atoms $B\neq 1,\Delta$ such that $C<B$ and the ascending link with the poset
  of atoms $B\neq 1,\Delta$ such that $C\not< B$. The former is a cone with
  cone-point $C$ (unless $C=\Delta$ which corresponds to $v=*$). To understand
  the latter, for each atom $D$ consider the subset ${\cal
    S}(D)=\cup\{B\mid B<D\}$ of the Coxeter sphere, where we now identify an atom
  $B$ with the intersection of the chamber $\pi(B)(Q)$ with the unit sphere in
  $V_\W$.  The complement of $1\cup \cup\{B\mid C<B\}$ is covered by the
  intersections with the sets ${\cal S}(D)$ for $C\not< D$. Proposition
  \ref{characterization} implies that the collection of ${\cal S}(D)$'s is
  closed under nonempty intersections. The elements are contractible (this can
  be seen by observing that for $C\not< D$ the convex ball ${\cal S}(D)$
  intersects $\cup\{B\mid C<B\}$ in its boundary (if at all), it contains the
  simplex 1, and the boundary of ${\cal S}(D)$ intersects 1 in a proper
  collection of faces, ie, a ball, whose complement in $\partial(1)$ is
  also a ball), and thus the poset of the cover (which can be identified with
  the ascending link) is homotopy equivalent to the underlying space
  $S^d\setminus (1\cup \cup\{B\mid C<B\})$. If the convex balls $1$ and $\cup\{B\mid C<B\}$
  are disjoint, then this complement is homotopy equivalent to $S^{d-1}$ and
  if they intersect then the complement is contractible. \end{proof}

\begin{remark} The proof also shows that the dualizing module
  $H_c^d(X(\G);\Z)$ is a free abelian group. In this situation,
$H_c^d(X(\G);\Z)$ is isomorphic to $$\bigoplus_{v\in V}
H^{d-1}(Lk_\uparrow(v,X(\G));\Z).$$ For details see
\cite{bf:duality}.\end{remark}

\begin{remark} \label{contractibility} 
There is another proof of the contractibility of $X(\G)$ that follows
  the lines of Charney's Fellow Traveller Property
  \cite{charney:biautomatic}. One shows that the function $\phi\co V\to
  V$ defined by $\phi(*)=*$ and $$\phi(A_1\cdot A_2\cdot \cdots \cdot
  A_{k-1}\cdot A_k)=A_1\cdot A_2\cdot\cdots\cdot A_{k-1}$$ extends to a
  simplicial map $X(\G)\to X(\G)$ which is homotopic to the
  identity. Then one observes that the image of any finite subcomplex
  of $X(\G)$ under a high iterate of $\phi$ is a point.
\end{remark}

\begin{remark} The connectivity at infinity of right-angled Artin
groups was studied by N. Brady and Meier \cite{bm:artin}.  These
groups act cocompactly on CAT(0) spaces and the distance from a point
is a good Morse function in that context. Brown and Meier
\cite{bm:duality} analyze which Artin groups of infinite type are
duality groups. \end{remark}

\subsection{NPC}

We wish to study the properties of the distance function on $X(\G)$, and
in particular its ``non-positively curved'' aspects. We start by
studying triangles formed by geodesics. It is convenient to
translate the triangle so that one vertex is $*$. Denote the other two
vertices by $v$ and $w$ (see figure 3).

\begin{prop} \label{down-and-up} All edges in the geodesic
from $v$ to $w$ whose orientation points towards $v$ occur at the
beginning of the path, and these are followed by the
edges oriented towards $w$.\end{prop}

In other words, the wordnorm first decreases and then increases along
the geodesic.

\begin{figure}[ht!]
\centerline{
\relabelbox\small
\epsfxsize0.4\textwidth\epsfbox{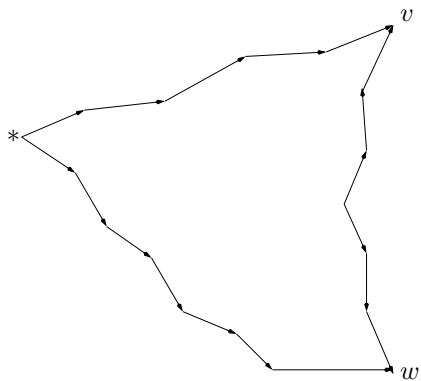}
\relabela <-3pt,0pt> {*}{$*$}
\relabel {v}{$v$}
\relabel {w}{$w$}
\endrelabelbox}
\caption{A geodesic triangle in $X(\G)$}
\end{figure}

\begin{proof} Suppose on the contrary that $x,y,z$ are (special
representatives of) consecutive vertices on the geodesic from $v$
to $w$ and that $||x||<||y||>||z||$.  Thus there are atoms $B$ and
$C$ such that $xB=y$ and $yC=z\Delta$. It follows that the first
atom in the normal form for $x(BC)=z\Delta=\Delta\overline z$ is
$\Delta$. Since $B\cdot C$ is a normal form, 
Proposition \ref{nf} implies that the first
atom in the normal form for $xB=y$ is $\Delta$, a
contradiction. \end{proof}

The following is the key ``non-positively curved'' feature of the
complex $X(\G)$. There are two versions, one for each metric. The
$d_{wd}$--version is more useful, since the inequality is strict. It
will allow us to run Cartan's fixed-point argument.

\begin{thm} \label{npc-D}
Suppose $p$ is a vertex on the geodesic
from $v$ to $w$ and $v\neq p\neq w$. Let $u$ be any vertex. Then
$$d_{at}(u,p)\leq\max\{d_{at}(u,v),d_{at}(u,w)$$ and
$$d_{wd}(u,p)<\max\{d_{wd}(u,v),d_{wd}(u,w)\}.$$ 
Moreover,
$$d_{wd}(u,p)\leq \max\{d_{wd}(u,v),d_{wd}(u,w)\}-
\min\{d_{at}(v,p),d_{at}(w,p)\}.$$
\end{thm}

\begin{proof} Without loss of generality $u=*$ and then the
proof follows from Proposition \ref{down-and-up}.\end{proof}

We remark that the function $d_{wd}(*,\cdot)$ is not necessarily convex
along $[v,w]$, cf $\A=\Z\times \Z\times\Z$.

\subsection{Centers}

Let $T$ be a finite set of vertices of $X(\G)$. The {\it circumscribed radius}
$r(T)$ is the smallest integer $r$ such that for some vertex $v\in X$,
$r=\max_{t\in T}d_{wd}(t,v)$. (Note the order of points $t$ and $v$ here.) A
vertex $v$ is said to be a {\it center} of $T$ if $r(T)=\max_{t\in
 T}d_{wd}(t,v)$. A set might have more than one center, but we have the next
best thing.

\begin{prop} If $v$ and $w$ are two centers of $T$, then $v$
and $w$ span an edge.\end{prop}

\begin{proof} Suppose the geodesic from $v$ to $w$ passes through
another vertex $z$. If $t\in T$, then
$$d_{wd}(t,z)<\max\{d_{wd}(t,v),d_{wd}(t,w)\}\leq r(T)$$ contradicting
the definition of
$r(T)$.\end{proof} 

We next have the Cartan fixed point theorem in our setting.

\begin{thm} \label{cartan}
Every finite subgroup $H<\G$ leaves a simplex of $X(\G)$
invariant (and fixes its barycenter). In particular, there are only
finitely many conjugacy classes of finite subgroups in $\G$.\end{thm}

\begin{proof} Let $T$ be an orbit of $H$. The set of centers of $T$ is
$H$--invariant and spans a simplex in $X(\G)$.\end{proof}

\subsection{Minsets of elements}

In the theory of isometric actions on nonpositively curved spaces an
important role is played by the minsets of isometries. These are the sets
of points that move the least under the given isometry. Here we
explore the analogous concept using the worddistance.

\begin{prop} Let $g\in \G$ and suppose $v,w$ are two vertices
of $X(\G)$ and $d_{wd}(v,g(v))\geq d_{wd}(w,g(w))$. Then for each vertex $z$
on the geodesic from $v$ to $w$ we have $d_{wd}(z,g(z))\leq
d_{wd}(v,g(v))$.\end{prop}

\proof We may assume that $v=*$, by left-translating (and
replacing $g$ by a conjugate). Let $w=X_1\cdot X_2\cdot \cdots\cdot
X_k$ be the normal form without $\Delta$'s. The element $g$ has the
form $g=B_1\cdot B_2\cdot \cdots\cdot B_m$ (no $\Delta$'s) or 
$g=B_1\cdot B_2\cdot \cdots\cdot B_m\Delta$. 

{\bf Case 1}\qua $g=B_1\cdot B_2\cdot \cdots\cdot B_m$. 

Note that the vertex $g(w)$ can be written as $B_1\cdots B_mX_1\cdots
X_k$, and also as $X_1\cdots X_ku$ where $||u||=d_{wd}(w,g(w))$. Since
we are assuming that $||u||\leq ||g||=d_{wd}(*,g(*))$, it follows that
we have equality in $\A^+$:
$$B_1\cdots B_mX_1\cdots X_k=X_1\cdots X_ku\Delta^r$$ for some $r\geq
0$. In particular, $X_1\cdots X_i<B_1\cdots B_mX_1\cdots X_i$ (by 
Proposition \ref{summary+}) and so we
have $$B_1\cdots B_mX_1\cdots X_i=X_1\cdots X_iu_i$$ for some $u_i\in
A^+$. In particular, $$d_{wd}(X_1\cdots X_i,g(X_1\cdots X_i))\leq
||u_i||=||B_1\cdots B_m||=||g||=d_{wd}(*,g(*)).$$ 

{\bf Case 2}\qua $g=B_1\cdot B_2\cdot \cdots\cdot B_m\Delta$.

The calculation is similar.  We now get $$B_1\cdots B_m\Delta X_1\cdots
X_k=X_1\cdots X_ku\Delta^r$$
and $r\geq 1$ since $||u||\leq ||g||=||B_1\cdots
B_m||$. Now push the lone $\Delta$ on the left-hand side all the way to the
right and cancel it to get $B_1\cdots B_m\overline X_1\cdots \overline
X_k=X_1\cdots X_ku\Delta^{r-1}$. Again by Proposition \ref{summary+} we get
$$X_1\cdots X_i<B_1\cdots B_m\overline X_1\cdots\overline X_i$$
so that
$$B_1\cdots B_m\overline X_1\cdots\overline X_i=X_1\cdots X_iu_i$$ for
some $u_i\in \A^+$ with $||u_i||=||B_1\cdots B_m||=||g||$. Passing again to
$\A/\<\Delta\>$ we obtain $g(X_1\cdots X_i)=X_1\cdots X_iu_i$ and so
$$d_{wd}(X_1\cdots X_i,g(X_1\cdots X_i))\leq ||u_i||=||g||.\eqno{\qed}$$

If $g\in \G$ and $d\in\Z$, let $T_d(g)=\{v\in X^{(0)}\mid d_{wd}(v,g(v))\leq
d\}$. Recall that $d_{wd}$ is not symmetric and note that writing
$d_{wd}(g(v),v)$ instead in the definition would give $T_d(g^{-1})$. The
minimal $d$ with $T_d(g)\neq \emptyset$ is the {\it displacement
length} of $g$, denoted $\lambda(g)$. The set $T_d(g)$ for this $d$ is
the {\it minset} of $g$, denoted $A_g$.

A set $T$ of vertices is {\it starlike} with respect to $t_0\in T$ if
for every $t\in T$ all vertices on the geodesic from $t_0$ to $t$
are in $T$. It can be shown, by the argument in Remark
\ref{contractibility}, that the span of a starlike set is
contractible. A set of vertices is {\it convex} if it is starlike with
respect to each of its vertices.

\begin{cor} The minset $A_g$ of each $g\in \G$ is convex. \qed\end{cor}

\begin{lemma}\label{displacement}
Let $g\in\A^+$ have normal form $B_1\cdot B_2\cdot \cdots\cdot
B_k$. Also suppose that the displacement of $*$ is not smaller than
the displacement of the vertex $v$ with special representative
$D_1\cdot D_2\cdot\cdots\cdot D_l$, ie, that $d_{wd}(*,g(*))\geq
d_{wd}(v,g(v))$. Then $D_1D_2\cdots D_l<B_1B_2\cdots B_kD_1D_2\cdots
D_l$. \end{lemma}

\begin{proof} Let $\tilde D=D_1D_2\cdots D_l$.
$d_{wd}(D_1\cdots
D_l,gD_1\cdots D_l)$ is computed as follows: let $\tilde E=E_1\cdot
E_2\cdot\cdots \cdot E_l$ be the normal form with $\tilde E\tilde
D=\Delta^l$, find the normal form of $\tilde E(B_1\cdots B_k\tilde
D)$, discard $\Delta$'s and count letters. Since $d_{wd}(*,g(*))\geq
d_{wd}(v,g(v))$, at least $l$ $\Delta$'s will have to be discarded,
ie $\Delta^l<\tilde EB_1B_2\cdots B_k\tilde D$. Write
$\Delta^l=\tilde E\tilde D$ and cancel $\tilde E$ to obtain $\tilde
D<B_1B_2\cdots B_k\tilde D$.
\end{proof}

\section{Applications}\label{sec:applications}

We now turn to the group-theoretic consequences of the discussion of
the geometric properties of $X(\G)$.

\subsection{Translation lengths}

The {\it translation length} of $g\in \G$, denoted $\tau(g)$, is
defined to be $$\tau(g)=\lim_{n\to\infty}\frac{||g^n(*)||}{n}.$$ See
\cite{gs:biautomatic} for basic properties of translation
functions (in the more general setting of finitely generated groups). 
An important feature of $\tau$ is that conjugate elements
have the same translation length and the restriction of $\tau$ to a
torsion-free abelian subgroup is a (non-symmetric) semi-norm. In
\cite{gs:biautomatic}, Gersten and Short prove that in biautomatic
groups elements of infinite order have nonzero translation lengths,
and this gives rise to many properties of such groups (eg a
biautomatic group cannot contain the Heisenberg group as a subgroup).
Finite type Artin groups are biautomatic by
\cite{charney:biautomatic}. Here we show that translation lengths of
torsion-free elements in $\G$ are bounded away from 0, and this gives
rise to additional group-theoretic properties.

Let $g\in \G$ be an element which is not conjugate to an atom, or an atom
followed by $\Delta$ (we shall se momentarily that this forces $g$ to have
infinite order). By $T_n$ denote the set $\{*,g(*),g^2(*),\cdots,
g^{n-1}(*)\}$. Let $r_n$ be the circumscribed radius of $T_n$. Since
$T_1\subset T_2\subset\cdots$ we have $r_1\leq r_2\leq
\cdots$. Suppose that $r_n=r_{n+1}$ for some $n$. Let $z$ be a center
for $T_{n+1}$. Then the ball of radius $r_n=r_{n+1}$ centered at $z$
covers both $T_n$ and $g(T_{n})$ and so $z$ is a center for both. But
then $g(z)$ is also a center for $g(T_n)$ and we conclude that $z$ and
$g(z)$ span an edge, so that $g$ is conjugate to an atom or an atom
followed by $\Delta$. Thus $r_n\geq n$ for all $n$ (and in particular $g$ has
infinite order). 

We now see that $\max_{i\leq n} d_{wd}(*,g^i(*))\geq n$, and in particular
there are infinitely many $n$ such that $d_{wd}(*,g^n(*))\geq n$. It
follows that the translation length $\tau(g)\geq 1$.

Recall that $\delta=length(\Delta)$.

\begin{thm} \label{trlengths}
The set of translation lengths of elements of $\G$
of infinite order is bounded below by $\frac{1}{2\delta}$.\end{thm}

\begin{proof} We have seen that if $g$ is of infinite order and not
conjugate to $B$ or $B\Delta$ for an atom $B$, then $\tau(g)\geq
1$. Thus if $g^{2\delta}$ is not conjugate to an atom or its inverse,
then $\tau(g)=\frac{1}{2\delta}\tau(g^{2\delta})\geq
\frac{1}{2\delta}$. 

The length homomorphism $length\co \A\to\Z$ descends to the homomorphism
$length\co \G\to\Z/2\delta$. Since conjugate elements have
the same length, $length(g^{2\delta})$ $=0\in\Z/2\delta$, and the only
atom with 0 length is the identity element, it follows that if
$g^{2\delta}$ is conjugate to an atom or its inverse, then $g$ has
finite order.
\end{proof}

\begin{cor} \label{fg}
Every abelian subgroup $A$ of $\G$ (or $\A$) is finitely
generated. \end{cor}

\begin{proof} Since $\G$ is virtually torsion-free, after passing to a
finite index subgroup of $A$ we may assume that $A$ is
torsion-free. Since the virtual cohomological dimension of $\G$ is
finite, it follows that $A$ is a subgroup of the finite-dimensional
vector space $W=A\otimes \Q$.  Since $\tau\co A\to [0,\infty)$ extends to
a (non-symmetric) norm on $W$ and its values on $A\setminus\{0\}$ are
bounded away from 0, it follows that $A$ is discrete in $W$. (Strictly
speaking, this is not a norm since $\tau$ is not symmetric. We could
symmetrize, or else work with a non-symmetric norm which is just as
good in the present context.)\end{proof}

\begin{cor} $\G$ does not have any infinitely divisible elements of
infinite order.\qed\end{cor}

The following corollary is a consequence of Theorem 3.5 of G. Conner's
work on translation lengths \cite{conner:trlengths}. For completeness,
we sketch the proof.

\begin{cor} Every solvable subgroup of $\G$ (or $\A$) is virtually
abelian.\end{cor}

\begin{proof} Let $H$ be a solvable subgroup of $\G$. Since $\G$ is
virtually torsion-free, we may assume that $H$ is torsion-free after
replacing it with a subgroup of finite index.  By induction on the
length of the derived series, the commutator subgroup $[H,H]$ is
virtually abelian, and hence finitely generated. Let $K\subset [H,H]$
be an abelian subgroup of finite index which is characteristic in
$[H,H]$. The conjugation action of $H$ on $K$ is virtually trivial,
since it preserves translation lengths (the image of the homomorphism
$H\to GL(K)$ consists of elements of finite order and is
therefore finite). Thus the action of $H$ on $[H,H]$ is virtually
trivial as well. After passing to a subgroup of $H$ of finite index,
we may assume that this action is trivial. The proof concludes by an
argument in \cite{gs:biautomatic}. We now claim that $[H,H]$
is the trivial group, and so (the new) $H$ is abelian. Indeed, suppose
that $x,y\in H$ and $z=[x,y]=xyx^{-1}y^{-1}\neq 1$. By our
construction, $z$ commutes with both $x$ and $y$, and it has infinite
order. We compute that $[x^n,y^n]=z^{n^2}$ and so the translation
length of $z$ is 0, a contradiction.
\end{proof}

\subsection{Finite subgroups of $\G$}

It is well-known to the experts that all finite subgroups of $\G$ are
cyclic, and in fact the kernel of the homomorphism $\G\to\Z/2\delta$
is torsion-free, where $\delta=||\Delta||$ and the homomorphism is the
length modulo $2\delta$. (Note that the argument of Theorem
\ref{trlengths} gives another proof of this fact.)

In this section we will use the geometric structure of $X(\G)$ to give
a classification of finite subgroups of $\G$ up to conjugacy.

\begin{thm} \label{cyclic}
Every finite subgroup $H<\G$ is cyclic. Moreover, after conjugation,
$H$ transitively permutes the vertices of a simplex $\sigma\subset
X(\G)$ that contains $*$ and $H$ has one of the following two forms:

{\bf Type 1}\qua The order of $H$ is even, say $2m$. It is generated by an atom
$B$. The vertices of $\sigma$ are $*,B,B^2,\cdots,B^{m-1}$ (all atoms)
and $B^m=\Delta$. Necessarily, $\overline B=B$ (since $\Delta$ fixes
the whole simplex). 

{\bf Type 2}\qua The order $m$ of $H$ is odd, the group is generated by
$B\Delta$ for an atom $B$, and the vertices $*,B,B\overline
B,B\overline BB,\cdots,(B\overline B)^{(m-1)/2}B$ (all atoms) are
permuted cyclically and faithfully by the group (so the dimension of
$\sigma$ is $m-1$). Since $m$ is odd, the square $B\overline B$ of the
generator also generates $H$. 
\end{thm}

An example of a type 1 group is $\<B\>$ for $B=\sigma_1\sigma_3\sigma_2$
in the braid group $B_4/\<\Delta^2\>$ (of order 4). An
example of a type 2 group is $\<\sigma_1\Delta\>$ in $B_3/\<\Delta^2\>$
(of order 3).

The key to this is:

\begin{lemma} \label{ordering}
The set of vertices of any simplex $\sigma$ in $X$
admits a cyclic order that is preserved by the stabilizer
$Stab(\sigma)<\G$. \end{lemma}

\begin{proof} We can translate $\sigma$ so that $*$ is one of its
vertices. Let the cyclic order be induced from the linear order 
$*<B_1<B_2<\cdots<B_k$ given
by the orientations of the edges of $\sigma$ (equivalently, by the
lengths of the atoms $B_i$). We need to argue that
the left translation by $B_i^{-1}$ is going to produce the same cyclic
order. We can write $B_{i+1}=B_iC_i$ for atoms $C_i$, so that the
vertices of $\sigma$ are $*,C_1,C_1C_2,\cdots,C_1\cdots C_k$. After
translation by $B_1^{-1}$ the vertices are
$B_1^{-1},*,C_2,C_2C_3,\cdots,C_2\cdots C_k$, where $B_1^{-1}$ should
be replaced by the special representative of the coset
$B_1^{-1}\<\Delta\>$. Let $Y$ be the atom with $C_1C_2\cdots
C_kY=\Delta$. Then $B_1^{-1}\Delta=C_1^{-1}\Delta=C_2C_3\cdots C_kY$
and $C_2C_3\cdots C_kY$ is this canonical representative (it is a
subword of $\Delta$, so it is an atom). Since $C_2\cdots C_k<C_2\cdots
C_kY$, it follows that the induced ordering on the vertices is a
cyclic permutation of the old one. Repeating this $i$ times gives the
claim. \end{proof}

\begin{proof}[Proof of Theorem \ref{cyclic}]
Recall that by Theorem \ref{cartan} there is an $H$--invariant
simplex $\sigma\subset X(\G)$.  By passing to a face if necessary and
conjugating, we may assume that $H$ acts transitively on the vertices
of $\sigma$, and that $*$ is a vertex of $\sigma$. Say there are $m$
vertices and choose $h\in H$ that rotates the simplex by one
unit. Note that $h^m$ is either $1$ or $\Delta$. Also, $h$ is either
an atom $B$ or else $h=B\Delta$ for an atom $B$. In the first case,
the vertices of the simplex are $*,B,BB,\cdots,B^{m-1}$ and these are
all atoms (since the arrow points from $B^i$ to $B^{i+1}$ for each
$i=0,1,\cdots,m-2$). Thus $1<B^{m-1}<\Delta$ and hence $B<B^m<B\Delta$,
so it follows that $B^m=\Delta$, and $H$ is of type 1. In the second
case, the vertices of $\sigma$ are $*,B,B\overline B,B\overline
BB,\cdots$ and they are similarly all atoms. If $m$ is even, the last
vertex in this sequence is $(B\overline B)^{m/2-1}B$, and this is
still an atom. It follows that $(B\overline B)^{m/2}=\Delta$. But then
$h^m=\Delta\in H$ fixes the whole simplex and in particular $B=\overline
B$ and $H$ is of type 1. So suppose $m$ is odd. Then the last vertex
of $\sigma$ is $(B\overline B)^{\frac{m-1}{2}}$ and we have
$(B\overline B)^{\frac{m-1}{2}}B=\Delta$. Multiplication by $\Delta$
on the right reveals that $h^m=1$ and $H$ is of type 2.
\end{proof}

\subsection{Normal abelian subgroups and centers of finite index subgroups}

We use minsets to prove:

\begin{thm} \label{infinite}
  Assume that the associated Coxeter group $\W$ is nonabelian and irreducibe.
  The action of $\G$ on itself by conjugation does not have nontrivial finite
  orbits. (Singletons consisting of central elements are trivial orbits.) In
  particular, the center of any finite index subgroup of $\G$ is either
  trivial or $\<\Delta\>$ (if the latter is central in $\G$).\end{thm}

We remark that $\Delta$ is central in $\G$ if and only if it is
central in $\A$ (if $\Delta g=g\Delta$ in $\G$ then {\it a priori} we
only get $\Delta g=g\Delta^{2m+1}$ in $\A$, but $m$ must be 0 by
length considerations) if and only if the bar
involution is trivial.

We now turn to the lemmas needed in the proof.

\begin{lemma} Suppose $A_1,\cdots,A_p$ is a finite collection of convex
sets in $V=\G/\<\Delta\>$ permuted by left translations. Then each $A_i$
is empty or all of $V$.\end{lemma}

\begin{proof} If not, we may assume $\emptyset\neq A_i\neq V$ for
each $i$, by discarding the copies of $\emptyset$ and $V$ from the
collection. Let $q\geq 1$ be the largest integer such that for some
$i_1<i_2<\cdots <i_q$ the intersection
$$\bigcap_{j=1}^q A_{i_j}\neq\emptyset.$$ Now pass to the collection
of $q$--fold intersections of the $A_i$'s. We can therefore assume that
the sets in the collection, still denoted $A_i$, are pairwise
disjoint. Since $\G$ acts transitively on $V$, the sets $A_i$
must cover $V$. Say $*\in A_1$. Let $\sigma$ be a generator of
$\A^+$ and consider the line $L=\{\sigma^j(*)\mid j\in\Z\}$. By
convexity, if $\sigma^j(*),\sigma^l(*)\in A_i$, then $\sigma^k(*)\in
A_i$ for $k$ between $j$ and $l$ (this is because the normal form of
$\sigma\sigma\cdots\sigma$ is
$\sigma\cdot\sigma\cdot\cdots\cdot\sigma$). We conclude that for $j$
large each $\sigma^j(*)$ belongs to the same $A_i$. Left translation
by $\sigma^{-j}$ takes $\sigma^j(*)\in A_i$ to $*\in A_1$ and it takes
$\sigma^{2j}(*)\in A_i$ to $\sigma^{j}(*)\in A_i$, and so
$i=1$. Convexity now implies that $\sigma^j(*)\in A_1$ for all $j\geq
0$, and in particular $\sigma(*)\in A_1$.

Let $v$ be now any vertex and assume $v\in A_1$. Choose $g\in \G$ with
$g(v)=*$. Thus $g(A_1)=A_1$ and we then have $g^{-1}(\sigma(*))\in
A_1$. The two choices of $g$ give $v\sigma\in A_1$ and
$v\overline\sigma\in A_1$. Since any vertex can be reached from $*$ by
successively right-multiplying by a generator, we conclude $A_1=V$, a
contradiction. 
\end{proof}

It is convenient to introduce the following notion.  Recall that
$x_1\cdot x_2\cdots x_k$ is a normal form if and only if
$x_i=\alpha(x_ix_{i+1})$ for all $i$. Motivated by this observation we
construct, following R Charney \cite{charney:biautomatic}, a finite
graph whose vertex set is the set of atoms not equal to 1 or $\Delta$,
and there is an arrow from $x$ to $y$ if $x=\alpha(xy)$. A nontrivial
normal form without $\Delta$'s is simply a finite directed path in
this graph. We refer to this graph as the {\it Charney graph}. So the
elements of $V$ are in 1--1 correspondence with oriented paths in the
Charney graph (with $*$ corresponding to the empty path, and atoms
corresponding to one point paths).

\begin{prop} \label{charney}
If $\W$ is irreducible, then any two atoms
$x,y$ not equal to 1 or $\Delta$ can be joined by an oriented path in
the Charney graph.
\end{prop}

\begin{proof}
As the first case, we assume that $x$ and $y$ are generators. Let
$x=g_1,g_2,\cdots g_k=y$ be a sequence of generators such that
successive elements do not commute (this sequence exists by the
irreducibility assumption). Then $$g_1\cdot g_1g_2\cdot g_2\cdot
g_2g_3\cdot g_3\cdots g_{k-1}\cdot g_{k-1}g_k\cdot g_k$$
is the desired path.

Next, we observe that if $x\neq 1$ is any atom, then there is a generator $g$
such that $x\cdot g$ is a normal form. Simply take $g$ to be the last
generator in a word representing $x$.

It remains to argue that if $y$ is any atom not equal to 1 or
$\Delta$, there is a generator $g$ and an oriented path from $g$ to
$y$. Let $S$ be the set of generators $a$ such that some word
representing $y$ begins with $a$. If $S$ consists of a single element
$a$, then $a\cdot g$ is a normal form and we are done.  We proceed by
induction on the cardinality of $S$. We first note that $S$ is a
proper subset of the generating set (if a chamber is separated from
the fundamental chamber $Q$ by every wall adjacent to $Q$, then the
chamber is the antipodal chamber $-Q$). Next, observe that
$\Delta(S)\cdot y$ is a normal form, where $\Delta(S)$ is the
$\Delta$--element in the subgroup generated by $S$. Let $c$ be a
generator not in $S$. Then $c\Delta(S)$ is an atom that can begin with
$c$ or with an element of $S$ that commutes with $c$ (this follows
from the cancellation law and the fact that if $y$ begins with $c$ and
$x$ then it begins with $\p(c,x,m(c,x))$). Since $c\Delta(S)\cdot
\Delta(S)\cdot y$ is a normal form, we have replaced the original
set $S$ by the set $S'$ consisting of $c$ and the elements of $S$ that
commute with $c$. 

{\bf Case 1}\qua $c$ can be chosen so that at least two elements
of $S$ don't commute with $c$. 

Then $card(S')<card(S)$ and we are
done by induction.

Recall that the Coxeter graph is a tree, and consider the forest
spanned by $S$. Two generators commute if and only if they are not
adjacent in the Coxeter graph.

{\bf Case 2}\qua The forest has more than one component.

If there are two components separated by a single vertex $c$, then clearly
$card(S')$ $<card(S)$ and we are done. Otherwise choose $c$ to be
adjacent to one component while separating it from another, and so
that the distance between the two components is as small as
possible. Then $S'$ has the same cardinality as $S$, but the
associated forest has two components that are closer together than in
the old forest. Repeating this procedure eventually produces two
components separated by a single vertex.

{\bf Case 3}\qua The forest is a tree and $card(S)>1$.

Then choose $c$ to be adjacent to a vertex in the tree. $S'$ has the
same cardinality as $S$ but the underlying forest has $>1$ component.
\end{proof}

\begin{lemma} \label{heck}
Suppose $\W$ is irreducible. Then for any chamber $R$
and any wall $W$ that abuts $R$ there is a normal form $D_1\cdot
D_2\cdot\cdots \cdot D_l$ such that
\begin{itemize}
\item $D_l$ is a single generator,
\item The gallery associated to $D_1D_2\cdots D_l$ that starts at the
fundamental chamber $Q$ ends at the chamber $R$ and last crosses wall
$W$.
\end{itemize}
Moreover, if $B\neq 1,\Delta$ is a given atom, the normal form can be
chosen so that $B\cdot D_1$ is a normal form.
\end{lemma}

\begin{proof}
Using the connectivity of the Charney graph, start with a normal form
$E_1\cdot E_2\cdot\cdots\cdot E_m$ so that $E_m$ is a single generator
(and $B\cdot E_1$ is a normal form).  Say $R'$ is the terminal chamber
of the gallery $E_1\cdots E_l$ and the last wall crossed is $W'$. If
$X$ is a generator that does not commute with $E_m$, then $E_1\cdot
E_2\cdot\cdots \cdot E_m\cdot E_m\cdot (E_mX)\cdot X$ is also a normal
form whose gallery ends at $R'$, but the last wall crossed is a
different wall from $W'$. By irreducibility, any two generators can be
connected by a sequence of generators with successive generators
noncommuting, and thus we can construct a gallery as above that ends
at $R'$ and last crosses any preassigned wall abutting $R'$. By repeating
the last atom in such a normal form, we can construct a similar gallery that
ends in any preassigned chamber adjacent to $R'$ and by iterating
these operations we can get to $R$ and $W$.
\end{proof}

\begin{lemma} If $\W$ is irreducible, then $A_g\neq V$ for
every nontrivial $g\in \G$, unless $g=\Delta$ and $\Delta$ is
central. \end{lemma}

We will actually prove a stronger form of this lemma. Consider the set
$H\subset \G$ of elements $g\in \G$ such that the function $v\mapsto
d_{wd}(v,g(v))$ defined on $V$ is bounded above. It is then clear that
$H$ is a normal subgroup of $\G$.

\begin{prop} \label{4.12} If $\W$ is irreducible, then $H$ is either
trivial or equal to $\<\Delta\>$ (if $\Delta$ is central).\end{prop}

\begin{proof} We will first argue that every element of $H$ has finite
order. This will imply that $H$ is finite (since $H\to\W$ will
then be injective).

Let $g\in H$ have infinite order. Notice that 
$||g^i(*)||\to\infty$ as $i\to\infty$. This is because
$g^i(*)=g^j(*)$ for $i\neq j$ would imply that $g^{i-j}$ equals 1 or
$\Delta$, so $g$ would have finite order. Fix a large number $N$ (to
be specified later) and replace $g$ by a power if necessary so that
$||g(*)||>N$. Further, by conjugating $g$ if necessary, we may assume
that $||g(*)||=d_{wd}(*,g(*))$ realizes the maximum of the (bounded) function
$v\mapsto d_{wd}(v,g(v))$.

As usual, we have either $g=B_1\cdot B_2\cdot\cdots\cdot B_k$ or
$g=B_1\cdot B_2\cdot\cdots\cdot B_k\Delta$. Note that $k>N/\delta$ is
large.

{\bf Case 1}\qua $g=B_1\cdot B_2\cdot\cdots\cdot B_k$.

If $\tilde D=D_1\cdot D_2\cdot \cdots\cdot D_l$ is a normal form such that
$B_k\cdot D_1$ is a normal form, then consider the vertex $D_1\cdots
D_l(*)=D_1\cdots D_l$. From Lemma \ref{displacement} we get $D_1D_2\cdots
D_l<B_1B_2\cdots B_kD_1D_2\cdots D_l$. Assume in addition that $l\leq k$.
Since the right-hand side is also a normal form, Lemma \ref{l1} implies that
$D_1D_2\cdots D_l<B_1B_2\cdots B_l$.  We will use this only when $l=k$.
Summarizing, if $B_k\cdot D_1$ is a normal form, then $D_1\cdots D_k<B_1\cdots
B_k$.  But we now argue that we can choose $\tilde D=D_1\cdot
D_2\cdot\cdots\cdot D_k$ so that $D_1$ follows $B_k$ (in the Charney graph)
but $\tilde D\not <B_1\cdots B_k$.

Say $M$ is an integer such that the atomnorm of the galleries
constructed in Lemma \ref{heck} is bounded by $M$.

Each atom $B\neq \Delta$ crosses at most $\delta-1$ walls in the
Coxeter sphere. Thus $B_1\cdots B_k$ crosses at most $k(\delta-1)$
walls. Therefore some wall $W$ is crossed by $B_1\cdots B_k$ at most
$k(\delta-1)/\delta$ times and the same is true for any initial piece
of $B_1\cdots B_k$. We will get a contradiction (for large $N$ and
$k$) by arguing that $\tilde D$ can be chosen so that it crosses a
preassigned wall at least $k-M$ times (contradiction arising if
$k-M>k(\delta-1)/\delta$ ie, if $k>\delta M$, and this can be
arranged if $N>\delta^2 M$).

Start with a normal form $D_1\cdot D_2\cdot\cdots\cdot D_l$ with
$l\leq M$ so that $B_k\cdot D_1$ is a normal form, $D_l$ is a single
generator and the last wall crossed is $W$. This is possible by Lemma
\ref{heck} (a left translation may be necessary before applying the
lemma since $B_1\cdots B_k$ may end at a chamber different from $Q$).
Then set $D_i=D_l$ for $i=l+1,l+2,\cdots, k$.

{\bf Case 2}\qua $g=B_1\cdots B_k\Delta$.

This is entirely analogous, except for some overlines, and is left to
the reader.

The proof of the proposition then follows from the following result.
\end{proof}

\begin{prop} If $\W$ is irreducible and $H<\G$ is a
  finite normal subgroup, then $H$ is trivial or $H=\<\Delta\>$ with
  $\Delta$ central.\end{prop}

\begin{proof} Since $H$ is finite, there is an $H$--invariant simplex in
  $X(\G)$. Since $H$ is normal, every vertex belongs to an
  $H$--invariant simplex. If $H=\<\Delta\>$, then $*$ is fixed by $H$,
  and hence every vertex is fixed by $H$. Thus the bar involution is
  trivial and $\Delta$ is central. If $H$ is nontrivial and not equal
  to $\<\Delta\>$, then $H$ must contain a nontrivial atom $B\neq\Delta$
  (the inverse of an element of the form $B\Delta$ is an atom). Left
  translation by $B$ moves every vertex to an adjacent vertex, and so
  $d_{wd}(v,Bv)<\delta$ for every vertex $v$. On the other hand, for
  any $g\in \G$ the numbers $d_{wd}(v,g(v))$ are all congruent $\mod
  \delta$ to each other, as $v$ ranges over the vertices of $X(\G)$,
  so in our situation we see that all displacements $d_{wd}(v,Bv)$ are
  equal to each other.  By the irreducibility of $\W$, there is a
  normal form $B\cdot B_1\cdot B_2\cdot\cdots\cdot B_k$ with $B_k$ any
  prechosen atom. Lemmas \ref{displacement} and \ref{l1} now imply
  that $B_1\cdots B_k<BB_1\cdots B_{k-1}$. In particular, we see that
  $||B_k||\leq ||B||$. By choosing $B_k$ to have length $\delta-1$, we
  conclude that $B$ has length $\delta-1$. There is then a generator
  $\sigma$ such that $\sigma B=\Delta$ and we conclude that
  $B^{-1}=\Delta\sigma\in H$. But $\Delta\sigma$ has infinite order,
  unless $\sigma$ and $\overline\sigma$ are the only generators.  (To
  see this, note that $X\cdot X\cdot X\cdot\cdots\cdot X$ is a normal
  form where $X=\p(\sigma,\overline\sigma,m(\sigma,\overline
  \sigma))$.)  If $\sigma=\overline\sigma$, then there is only one
  generator and $G=\<\Delta\>\cong \Z_2$ so there is nothing to prove.
  If $\sigma\neq\overline\sigma$, then each half of the Artin relation
  $\sigma\overline\sigma\sigma\cdots=
  \overline\sigma\sigma\overline\sigma\cdots$ has length equal to an
  odd integer $m\geq 3$. In particular, $\sigma\overline\sigma\in H$
  is an atom, so it must have length $\delta-1$, and this forces $m=3$.
  We are now reduced to the classical braid group on 3 strands modulo
  the center. To finish the argument, note that in this case
  $\sigma\overline\sigma\cdot \overline\sigma$ is a normal form, while
  $\overline\sigma\not <\sigma\overline\sigma$.
\end{proof}

\begin{proof}[Proof of Theorem \ref{infinite}] Let $g_1,\cdots,g_n$ be a finite orbit
under conjugation. Note that there is $K>0$ such that $d_{wd}(v,g_i(v))\leq K$
for all $i=1,2,\cdots,n$ and all vertices $v$. Indeed, we can take $K=\max
d_{wd}(*,g_i(*))$, for then 
$$d_{wd}(v,g_i(v))=d_{wd}(h(*),g_ih(*))=d_{wd}(*,h^{-1}g_ih(*))=
d_{wd}(*,g_j(*))$$
where $h$ is chosen so that $h(*)=v$. It now follows from Proposition
\ref{4.12} that
each $g_i$ is central.\end{proof}

The following corollary answers a question of Jim Carlson. It
motivated the construction of $X(\G)$ and the analysis of its
geometric properties.

\begin{cor} Assume that the associated Coxeter group $\W$ is
  irreducible. Let $A$ be a normal abelian subgroup
  of $\G$. Then $A$ is trivial or $\<\Delta\>$ (and in the latter case
  $\Delta$ is central).
\end{cor}

\begin{proof} $\G$ acts on $A$ by conjugation. If $A$ is not as in the
  conclusion, then this action has infinite orbits
  by Theorem \ref{infinite}. The abelian group
  $A$ is finitely generated by Corollary \ref{fg}. The translation
  length function induces a norm on the free abelian group
  $A/\text{torsion}$. The induced action of $\G$ must preserve this
  norm and it still has infinite orbits, a contradiction.
\end{proof}

\section{The space at infinity}

We now construct a ``space at infinity'' of $X(\G)$ and examine the
basic properties of the action of $\G$. In a joint work with Mark
Feighn it will be shown that elements of infinite order have periodic
points at infinity. This is to be regarded as the analog of the space
of projectivized geodesic measured laminations for the case of mapping
class groups (when $\A$ is a braid group, $\G=\A/\Delta^2$
is a mapping class group).

\subsection{Definition}

Recall that a normal form without $\Delta$'s is an oriented path in
the Charney graph.

\begin{definition} An {\it admissible itinerary} is an infinite
directed path $X=x_1\cdot x_2\cdot x_3\cdots$ in the Charney graph.
\end{definition}

We denote by $\Omega$ the set of all admissible itineraries and topologize it
in the usual fashion: two are close if they agree for a long time. Thus
$\Omega$ is a totally disconnected compact metrizable space, and it is
nonempty unless $\A$ is trivial or $\Z$. If $\W$ is irreducible and
nonabelian, then $\Omega$ is a Cantor set, by Proposition \ref{charney} plus
the observation that there are vertices with at least two outgoing edges (eg
if $a,b\in S$ don't commute then there are oriented edges from $a$ to both $a$
and $ab$). We now describe an action of $\A^+$ on $\Omega$. Let $X\in\Omega$
be as above and let $g\in \A^+$. Observe that by Proposition \ref{nf}
$y_1=\alpha(gx_1)=\alpha(gx_1x_2\cdots x_k)$ for any $k\geq 1$.  Similarly, by
Proposition \ref{summary+} the second atom $y_2$ in the normal form for
$gx_1x_2$ equals the second atom in the normal form for $gx_1x_2\cdots x_k$
for any $k\geq 2$.  Continuing in this fashion, we see that the normal forms
for $gx_1x_2\cdots x_k$ ``converge'' as $k\to\infty$ to an infinite sequence
$y_1\cdot y_2\cdot y_3\cdots$ so that any finite initial piece is a normal
form. This may not be an admissible itinerary since $y_i$ might be $\Delta$.
However, observe that only finitely many $y_i$'s can be $\Delta$. They all
occur at the beginning and their number is no larger than the atomnorm of $g$.
If the number of $\Delta$'s is even we erase them, and if the number is odd we
erase them and replace all the remaining $y_i$ by their ``conjugate''
$\overline y_i$. Intuitively, we think of pushing all the $\Delta$'s off to
infinity, much in the same way as we calculate the special representative of a
vertex of $X(\G)$.

Since $\Delta^2$ acts trivially, there is an induced action of $\G$ on
$\Omega$. The action is continuous (this follows from Lemma
\ref{summary+}).

\subsection{Faithfulness and minimality}

\begin{prop} If $\A$ is not $1$, $\Z$, or $\Z\times
\Z$, then the kernel of the action of $\G=\A/\Delta^2$ 
on $\Omega$ is either trivial or $\<\Delta\>$ (if $\Delta$ is central).\end{prop}

\begin{proof} If $\Delta$ is not central, then there is a generator
$a$ such that $\overline a=b\neq a$. Then $a^\infty:=a\cdot a\cdot
a\cdots\in\Omega$ is not fixed by $\Delta$; indeed its image is
$b^{\infty}$.

Let $g\in \A^+$ be nontrivial. Assume first that the normal form of
$g=x_1\cdot x_2\cdots x_k$ does not have any $\Delta$'s. Let $a$ be a
generator such that $x_k\cdot a$ is a normal form (ie, so that $x_ka$
is not an atom --- eg $a$ can be the last letter in a word
representing $x_k$). Then $a^\infty:=a\cdot a\cdot a\cdots\in\Omega$
is fixed by $g$ only if $g=a^k$ is a power of $a$. If $b$ is a
generator that does not commute with $a$, then $a^k(b^\infty)=a\cdot
a\cdots a\cdot ab\cdot b\cdot b\cdots$ and $b^\infty$ is not fixed by
$g$. Finally, if $b$ is a generator distinct from $a$ and commuting
with $a$, but there are other generators (so that $\Delta\neq ab$),
then $a^k(b^\infty)=ab\cdot ab\cdots ab\cdot b\cdot b\cdots$, so
$b^\infty$ is not fixed in that case either.

Finally, assume that the normal form of $g$ is $\Delta\cdot
x_1\cdot\cdots\cdot x_k$ with $k\geq 1$ and $x_1\neq \Delta$. Again
take $a$ to be a generator with $x_k\cdot a$ a normal form. If $g$
fixes $a^\infty$, then $\overline x_1\cdot \cdots\cdot \overline
x_k\cdot\overline a\cdot\overline a\cdot\cdots=a\cdot a \cdot\cdots$
and so $g=\Delta a^k$ and $\overline a=a$. If $b$ is a generator such
that $\overline b$ does not commute with $a$, then $\Delta a^k
(b^\infty)=a^k\overline b^\infty=a\cdot a\cdot\cdots\cdot a\cdot
a\overline b\cdot \overline b\cdot\cdots$ and so $b^\infty$ is not
fixed. Finally, if $\A$ is abelian but has more than 2 generators and
if $b$ is a generator distinct from $a$, then $\Delta a^k
b^\infty=ab\cdot ab\cdot\cdots\cdot ab \cdot b\cdot b\cdot\cdots$.
\end{proof}

\begin{prop} Suppose that $\W$ is irreducible 
and nonabelian. Then each
orbit in $\Omega$ is dense. In particular, there is no proper closed
invariant subset of $\Omega$.
\end{prop}

\begin{proof} Let $X=x_1\cdot x_2\cdots$ and $Y=y_1\cdot y_2\cdots$ be two
points in $\Omega$. We need to construct some $g\in \A^+$ such that
$gX$ and $Y$ agree in the first $n$ slots. By Proposition
\ref{charney} there is a finite directed path $y_n\cdot z_1\cdots
z_k\cdot x_1$ from $y_n$ to $x_1$. Take $g=y_1y_2\cdots y_n z_1\cdots z_k$.
\end{proof}


\begin{thebibliography}

\bibitem{bb:morse} {\bf M~Bestvina}, {\bf N~Brady}, \emph{Morse theory
and finiteness properties of groups}, Invent. Math. {109} (1997)
445--470

\bibitem{bf:duality} {\bf M~Bestvina}, {\bf M~Feighn}, \emph{The
topology at infinity of ${O}ut({F}_n)$}, preprint (1997)

\bibitem{bh:cat0} {\bf M~Bridson}, {\bf A~Haefliger}, \emph{Metric
spaces of non-positive curvature}, manuscript of a book, in progress

\bibitem{bm:artin} {\bf N~Brady}, {\bf J~Meier}, \emph{Connectivity at
infinity for right-angled {A}rtin groups}, preprint

\bibitem{bm:duality} {\bf K\,S Brown}, {\bf J~Meier}, \emph{Improper
actions and higher connectivity at infinity}, in progress

\bibitem{bourbaki:coxeter} {\bf N~Bourbaki}, \emph{{G}roupes et
{A}lgebres de {L}ie}, Masson, Paris (1981) Chapters {IV--VI}

\bibitem{t.brady:artin} {\bf T~Brady}, \emph{Artin groups of finite
type with three generators}, preprint

\bibitem{brieskorn:bourbaki} {\bf E~Brieskorn}, \emph{Sur les groupes
de tresses}, S\' emininaire Bourbaki 24e ann\' ee, 1971/72, no. 401,
LNM, vol. 317, Springer--Verlag (1973) 21--44

\bibitem{bs:fartin} {\bf E~Brieskorn}, {\bf K~Saito}, \emph{Artin
{G}ruppen und {C}oxeter {G}ruppen}, Invent. Math. {17} (1972)
245--271

\bibitem{cd:browder} {\bf R~Charney}, {\bf M~Davis}, \emph{Finite
${K}(\pi,1)$s for {A}rtin groups}, from: ``Prospects in {T}opology'',
Proceedings of a Conference in honor of William Browder, (F~Quinn,
editor), Annals of Math. Studies, vol. 138, Princeton University Press
(1995) 110--124.

\bibitem{cd:jams} {\bf R~Charney}, {\bf M~Davis}, \emph{The
${K}(\pi,1)$--problem for hyperplane complements associated to
infinite reflection groups}, Jour. Amer. Math. Soc.  {8} (1995)
597--627

\bibitem{charney:biautomatic} {\bf R~Charney}, \emph{Artin groups of
finite type are biautomatic}, Math. Ann.  {292} (1992) 671--683

\bibitem{charney:symmetric} {\bf R~Charney}, \emph{Geodesic automation
and growth functions for {A}rtin groups of finite type},
Math. Ann. {301} (1995) 307--324

\bibitem{conner:trlengths} {\bf G~Conner}, \emph{Discreteness
properties of translation numbers in solvable groups}, preprint

\bibitem{deligne:fartin} {\bf P~Deligne}, \emph{Les immeubles des
groupes de tresses g\'en\'eralis\'es}, Invent. Math. {17} (1972)
273--302

\bibitem{garside:braids} {\bf F\,A Garside}, \emph{The braid groups
and other groups}, Quart. J. Math. Oxford {20} (1969) 235--254

\bibitem{gs:biautomatic} {\bf S~Gersten}, {\bf H~Short},
\emph{Rational subgroups of biautomatic groups}, Ann.  Math. {105}
(1991) 641--662

\bibitem{squier:artin} {\bf C\,C Squier}, \emph{The homological
algebra of {A}rtin groups}, Math. Scand.  {75} (1994) 5--43

\end{thebibliography}
\end{document}